\documentstyle[12pt,amsfonts]{article}

\begin{document}
\title{On Solutions of Three \\Quasi-geostrophic Models}

\author{ 
 Jiahong Wu\thanks{Partially supported by NSF grant DMS 9971926.} 
\\Department of Mathematics\\
 Oklahoma State University \\
 }
\date{}
\maketitle

\newtheorem{thm}{Theorem}[section]
\newtheorem{cor}[thm]{Corollary}
\newtheorem{prop}[thm]{Proposition}
\newtheorem{define}[thm]{Definition}
\newtheorem{rem}[thm]{Remark}
\newtheorem{example}[thm]{Example}
\newtheorem{lemma}[thm]{Lemma}
\def\theequation{\thesection.\arabic{equation}}

\vspace{.25in}
\noindent {\bf Abstract}. We consider the quasi-geostrophic model and its 
two different regularizations. Global regularity results are established
for the regularized models with critical or sub-critical indices. 
The proof (\cite{Ey0},\cite{CET}) of Onsager's conjecture \cite{On}
concerning weak solutions of the 3D Euler equations
and the notion of
dissipative solutions of Duchon and Robert \cite{DR} are extended to 
weak solutions of the quasi-geostrophic equation.

\vspace{.3in}
\noindent AMS (MOS) Numbers: 86A05, 35K55, 35Q35, 76U05

\vspace{.2in}
Regularized
quasi-geostrophic equation, Weak solutions, Global regularity.

\newpage

\section{Introduction}
\setcounter{equation}{0}
\label{sec:1}

Consider the two dimensional (2D) quasi-geostrophic (QG) equation
\begin{equation}\label{model1}
\theta_t + u\cdot \nabla \theta =0
\end{equation}
and its two different regularizations
\begin{equation}\label{model2}
\theta_t  + u\cdot \nabla \theta 
+ \kappa (-\Delta)^{\alpha} \theta=0
\end{equation}
and 
\begin{equation}\label{model3}
\theta_t + u\cdot \nabla \theta 
+ \mu (-\Delta)^{\alpha} \theta_t =0,
\end{equation}
where $\theta(x,t)$ is a real-valued function of $x$ and $t$,
$0\le \alpha\le 1$, $\kappa>0$ and $\mu>0$ are real numbers.
The advective velocity $u$ in these equations is determined from $\theta$ by  
a stream function $\psi$ via the auxiliary relations 
\begin{equation}\label{u}
u= (u_1,u_2) =\left(-\frac{\partial \psi}{\partial x_2},
\frac{\partial \psi}{\partial 
x_1}\right)\qquad \mbox{and}\qquad (-\Delta)^\frac12 \psi=-\theta. 
\end{equation}
Interest will mainly focus on the  behavior of solutions 
of the initial value problems (IVP)  
for these equations wherein 
\begin{equation}\label{init}
\theta(x,0) =\theta_0(x), \qquad\mbox{is specified}.
\end{equation}
To avoid questions regarding boundaries, we will assume periodic 
boundary conditions with period box $\Omega=[0,2\pi]^2$. 

\vspace{.15in}
Equations (\ref{model1}) and (\ref{model2}) are special cases of 
the general quasi-geostrophic 
approximations \cite{Pe} for atmospheric and oceanic fluid flow 
with small Rossby and Ekman numbers.
The variable $
\theta$ represents potential temperature, $u$ is the fluid velocity 
and $\psi$ can be identified with the pressure. 
Equation (\ref{model1}) is an important example 
of a 2D active scalar with a specific structure most
closely related to the 3D Euler equations while 
the equation in (\ref{model2}) with $\alpha=\frac12$ is the 
dimensionally correct analogue of the 3D Navier-Stokes
equations. These equations have recently been intensively investigated 
because of both 
their mathematical importance and their potential for applications in 
meteorology and oceanography (\cite{Pe},\cite{HPGS},\cite{CMT},\cite{Re},
\cite{OY}, \cite{CNS}, \cite{Co},\cite{CW1}).

\vspace{.15in} 
In modeling long waves in nonlinear dispersive media, Benjamin,Bona 
and Mahony \cite{BBM} introduced the BBM equation 
$$
u_t +u_x +uu_x - u_{xxt} =0
$$
as an alternative to the KdV equation 
$$
u_t +u_x +uu_x +u_{xxx} =0.
$$
Equation (\ref{model3}) to (\ref{model2}) is like the BBM to the KdV 
equation and our motivation for proposing such a model for study comes 
from the effects of regularizations on the global 
regularity of weak solutions and the potential applications of this new model 
in geophysics. 

\vspace{.15in}
These QG models appear to be simpler than the 
3D hydrodynamics equations, but contain many of 
their difficult features. For instance, solutions of 
(\ref{model1}) and (\ref{model2}) exhibit strong nonlinear 
behavior, strikingly analogous to that of the potentially singular 
solutions of the 3D hydrodynamics equations \cite{CMT}. 
Although progress has been made in the past several years (
\cite{CMT},\cite{Re},
\cite{OY}, \cite{CNS}, \cite{Co},\cite{CW1}),
the theory remains fundamentally incomplete. In particular, it is not 
known whether or not weak solutions of (\ref{model2}) are regular for 
all time when $\alpha$ is equal to the critical index $\frac12$. 
The critical case regularity issue turns out to be extremely difficult and 
was labeled by S. Klainerman \cite{Kl}
 as one of the most challenging PDE problems 
of the 21st Century. In Section \ref{sec:2} we explore how far one 
can go toward a regularity proof in the critical case 
and what are the weakest assumptions  needed to fill the gap. 

\vspace{.15in}
In Section \ref{sec:3} we solve the global regularity problem  
for equation (\ref{model3}) with $\alpha\ge \frac12$. We first construct a 
local solution $\theta$ in $H^s$ with $s>1$ and then derive explicit 
bounds on the norms of all derivatives of the solution. 
From this we infer that for 
$\alpha>\frac12$ the local solution $\theta$ remains bounded in $H^s$ 
for all time and thus no finite-time singularity can occur in this case.  
For the critical index $\alpha=\frac12$, global smoothness results
are established under  
assumptions that are much weaker than those needed to guarantee 
regularity for equation (\ref{model2}). This leads us to conclude that 
solutions of (\ref{model3}) are better behaved and thus (\ref{model3}) 
constitutes a reasonable alternative to (\ref{model2}).

\vspace{.15in}
As is well-known, weak solutions of the 3D hydrodynamics equations in
general only satisfy an energy inequality rather than equality. But 
Onsager conjectured in \cite{On} that weak solutions of the 3D Euler equations 
in a H\"{o}lder space $C^\gamma$ with exponent 
$\gamma>\frac13$ should conserve 
energy. In \cite{Ey0} Eyink proved energy conservation for 
weak solutions in a strong form of H\"{o}lder space 
$C^\gamma_*$ ($\gamma>\frac13$), in which the norm is defined in terms of 
absolute Fourier coefficients. Constantin, E and 
Titi provided a proof for the sharp version of Onsager's 
conjecture in \cite{CET}.  
Section \ref{sec:4} is concerned with solutions of the QG 
equation (\ref{model1}). First we verify Onsager's conjecture for weak 
solutions of the QG equation, extending the result of Constantin, E and Titi.
Then the QG equation is shown to possess the dissipative weak solutions, a 
notion proposed by Duchon and Robert \cite{DR}. Finally the two models 
(\ref{model1}) and (\ref{model3}) are proven to be close by considering
the limit of (\ref{model3}) as $\mu\to 0$. This provides 
further evidence 
for the validity of (\ref{model3}). 

\vspace{.15in}
We now review the notations used throughout the sequel. 
The Fourier transform $\widehat{f}$ of a tempered distribution $f(x)$ on 
$\Omega$ is defined as 
$$
\widehat{f}(k) =\frac1{(2\pi)^2 }\int_{\Omega} f(x) e^{-ik\cdot x} dx.
$$
We will denote the square root of the Laplacian $(-\Delta)^\frac12$
by $\Lambda$ and obviously 
$$
\widehat{\Lambda f}(k) = |k| \widehat{f}(k).
$$
More generally, $\Lambda^\beta f$ for $\beta\in {\Bbb R}$ 
can be identified with the Fourier series
$$
\sum_{k\in{\Bbb Z}^2} |k|^\beta \widehat{f}(k) e^{ik\cdot x}. 
$$
The equality relating $u$ to $\theta$ in (\ref{u}) 
can be rewritten in terms of 
periodic Riesz transforms
$$
u=\left(\partial_{x_2} \Lambda^{-1} \theta,\, -\partial_{x_1} 
\Lambda^{-1} \theta\right) = (-{\cal R}_2 \theta,\, {\cal R}_1\theta),
$$
where ${\cal R}_j, \,j=1,2$ denotes the Riesz transforms defined by
$$
\widehat{{\cal R}_j f}(k) = -i\frac{k_j}{|k|}\widehat{f}(k), \quad
k\in {\Bbb Z}^2\setminus \{0\}.
$$

\vspace{.15in}
$L^p(\Omega)$ denotes the space of the $p$th-power integrable functions
normed by
$$
|f|_{p} =\left(\int_{\Omega} |f(x)|^p dx \right)^\frac1p.
$$
For any tempered distribution $f$ on $\Omega$ and $s\in {\Bbb R}$, we
define 
$$
\|f\|_s = |\Lambda^s f|_2\,
=\left(\sum_{k\in {\Bbb Z}^2} |k|^{2s} |\widehat{f}(k)|^2
\right)^\frac12
$$ 
and $H^s$ denotes the Sobolev space of all $f$ for which $\|f\|_s$ is 
finite. For $1\le p\le \infty$ and $s\in {\Bbb R}$, the 
space $\frak{L}^p_s(\Omega)$ is a subspace of $L^p(\Omega)$, consisting
of all $f$ which can be written in the form $f=\Lambda^{-s} g$, $g\in L^p
(\Omega)$ and the ${\frak L}^p_s$ norm of $f$ is defined to be the $L^p$
norm of $g$, i.e.,
$$
\|f\|_{p,s} = |g|_p.
$$

\vspace{.2in}
\section{Dissipative QG Equation}
\setcounter{equation}{0}
\label{sec:2}

In this section we focus on the nonlinear behavior of weak solutions of the 
initial-value problem (IVP) for the dissipative QG equation
\begin{equation}\label{dqg}
\left\{
\begin{array}{ll}
\theta_t  + u\cdot \nabla \theta 
+ \kappa (-\Delta)^{\alpha} \theta=f,\quad 
& (x,t)\in \Omega\times [0,\infty),\\\cr
u=(u_1,u_2) = (-{\cal R}_2 \theta,\, {\cal R}_1\theta),\quad
& (x,t)\in \Omega\times [0,\infty),\\\cr
\theta(x,0) =\theta_0(x),\quad
& x\in \Omega,\cr
\end{array}
\right.
\end{equation}
where $0\le \alpha\le 1$ and $\kappa>0$ are real numbers. We establish 
nonlinear estimates which characterize the regularity of 
weak solutions of the IVP (\ref{dqg}) 
with $\alpha$ greater than or equal to the critical index $\frac12$.

\vspace{.15in}
For $0\le \alpha\le 1$, weak solutions of the IVP (\ref{dqg}) are known 
to exist globally in time \cite{Re}. More precisely, 
for any $T>0$,\, $\theta_0\in L^2$ and $f\in L^1([0,T];L^2)$, there exists a 
weak solution $\theta\in L^\infty([0,T]; L^2)\cap L^2([0,T];H^\alpha)$
satisfying 
$$
|\theta(\cdot,t)|_2^2 + \int_0^t \|\theta(\cdot,\tau)\|_{\alpha}^2 d\tau 
\le \left[|\theta_0|_2 + C \int_0^t |f(\cdot,\tau)|_2 d\tau\right]^2
$$ 
Furthermore, if $\theta_0 \in L^q$ anf $f\in L^1([0,T];L^q)$ 
for $1<q\le \infty$, then the  maximum principle
$$
|\theta(\cdot,t)|_q \le |\theta_0|_q + \int_0^t |f(\cdot,\tau)|_q d\tau
$$
holds for any $t\le T$.

\vspace{.15in}
Let $s\ge 0$. We now estimate $\|\theta\|_s=|\Lambda^s \theta|_2$.
If we take the inner product of $\Lambda^{2s}\theta$ with the first equation 
in (\ref{dqg}), we obtain 
\begin{equation}\label{root}
\frac12 \frac{d}{dt} |\Lambda^s \theta|^2_2 + \kappa |\Lambda^{s+\alpha}
\theta|^2_2  = (\Lambda^{2s}\theta,f) - (\Lambda^{2s}\theta, u\cdot
\nabla \theta) 
\end{equation}
The first term on the right hand side is bounded above by
\begin{equation}\label{easy}
|(\Lambda^{2s}\theta,f)|
\le |\Lambda^{s+\alpha}\theta|_2 | \Lambda^{s-\alpha}f|_2 
\le \frac{\kappa}{2}|\Lambda^{s+\alpha}\theta|^2_2 
    + \frac{1}{2 \kappa}|\Lambda^{s-\alpha}f|^2_2
\end{equation}
For the second term, we have 
\begin{equation}\label{hard}
|(\Lambda^{2s}\theta, u\cdot\nabla \theta)| 
=|(\Lambda^{2s}\theta, \nabla (u\theta)|
\le |\Lambda^{s+\beta}\theta|_2 \,|\Lambda^{s+1-\beta}(u\theta)|_2
\end{equation}
where $\beta\le \alpha$ remains to be determined.
To proceed, we need the calculus inequality 
\begin{equation}\label{cal}
|\Lambda^\gamma(FG)|_r 
\le C \Big[|\Lambda^\gamma F|_p |G|_q + |F|_q |\Lambda^\gamma G|_p\Big],
\end{equation}
where $\gamma>0$, $1<r\le p \le \infty$ and $1/r=1/p+1/q$. The use of
(\ref{cal}) gives
\begin{equation}\label{mid}
|(\Lambda^{2s}\theta, u\cdot\nabla \theta)|
\le C |\Lambda^{s+\beta}\theta|_2\left[|\Lambda^{s+1-\beta}u|_p  
|\theta|_q + |\Lambda^{s+1-\beta}\theta|_p |u|_q\right],
\end{equation}
where $2<q\le \infty$ and $1/p+1/q=1/2$. Taking into account of the second
equation in (\ref{dqg}) and the inclusion $H^{s+2-\frac2p-\beta}\subset 
{\frak L}^p_{s+1-\beta}$, we have 
$$
|\Lambda^{s+1-\beta}u|_p \le |\Lambda^{s+1-\beta}\theta|_p
\le |\Lambda^{s+2-\frac2p-\beta}\theta|_2.
$$
It then follows from (\ref{mid}) that 
$$
|(\Lambda^{2s}\theta, u\cdot\nabla \theta)|
\le C (|\theta|_q + |u|_q) |\Lambda^{s+\beta}\theta|_2
|\Lambda^{s+2-\frac2p-\beta}\theta|_2. 
$$
In the above, $\beta$ is essentially arbitrary and we may choose  
$$
\beta =1-\frac1p=\frac12+\frac1q
\qquad\mbox{so that $s+\beta=s+2-\frac2p-\beta$}.
$$
Hence 
\begin{equation}\label{last}
|(\Lambda^{2s}\theta, u\cdot\nabla \theta)|
\le C (|\theta|_q + |u|_q) |\Lambda^{s+\beta}\theta|_2^2.
\end{equation}

\noindent Combining (\ref{root}),(\ref{easy}) and (\ref{last}), we have 
\begin{equation}\label{key}
\frac{d}{dt} |\Lambda^s \theta|^2_2 + \kappa |\Lambda^{s+\alpha}
\theta|^2_2  \le \frac{1}{\kappa}|\Lambda^{s-\alpha}f|^2_2 
+ C_0(|\theta|_q + |u|_q) |\Lambda^{s+\beta}\theta|_2^2
\end{equation}
for any $s\ge 0$, $2<q\le\infty$ and $\beta=1/2+1/q\le \alpha$.

\vspace{.15in}
We now prove that weak solutions of the IVP (\ref{dqg}) with $\alpha>\frac12$
are actually regular. More precisely, we have the following theorem.

\begin{thm}\label{big}
Let $\alpha>\frac12$ and  
assume for $T>0$, $s>0$ and $2<q<
\frac{2}{2\alpha-1}$ 
$$\theta_0\in H^s\cap L^q\quad\mbox{and}\quad f\in L^1([0,T];L^2)\cap 
L^1([0,T];L^q)\cap L^2([0,T]; H^{s-\alpha}). $$
Then any weak solution $\theta$ of the IVP (\ref{dqg}) are regular in the
sense that  
$$
\theta \in L^\infty([0,T]; H^s) \cap L^2([0,T];H^{s+\alpha}) 
$$
\end{thm}
{\it Proof}.\quad The idea of the proof is to obtain from 
(\ref{key}) a closed differential
equality for $|\Lambda^s\theta|^2$. Since $u$ is essentially the Riesz
transform of $\theta$ from the second equation in (\ref{dqg}), we have
for $1<q<\infty$
\begin{equation}\label{ugood}
|u(\cdot,t)|_q \le |\theta(\cdot,t)|_q 
\le |\theta_0|_q + \int_0^t|f(\cdot,\tau)|_q d\tau.
\end{equation}
To eliminate the occurrence of $|\Lambda^{s+\beta}\theta|_2$, we use one
of the inequalities of Gagliardo and Nirenberg
$$
|\Lambda^{s+\beta}\theta|_2 \le C |\Lambda^{s+\alpha}
\theta |_2^{\frac\beta\alpha}
\,\,|\Lambda^s\theta|^{1-\frac\beta\alpha}_2. 
$$
We use H\"{o}lder's inequality to find 
\begin{equation}\label{bgood}
|\Lambda^{s+\beta}\theta|_2^2 \le \frac{\kappa}{2}
|\Lambda^{s+\alpha}\theta |_2^2 + \frac{C}{\kappa}|\Lambda^s\theta|^2_2.
\end{equation}
Inserting (\ref{ugood}) and (\ref{bgood}) into (\ref{key}) 
and canceling a factor of
$|\Lambda^{s+\alpha}\theta |_2^2$,
\begin{equation}\label{final}
\frac{d}{dt} |\Lambda^s \theta|^2_2 + \frac{\kappa}{2}|\Lambda^{s+\alpha}
\theta|^2_2  \le \frac{1}{\kappa}|\Lambda^{s-\alpha}f|^2_2 
+ \frac{C}{\kappa}|\Lambda^s\theta|^2_2,
\end{equation}
where $C$ only depends on $|\theta_0|_q$ and 
$\int_0^T|f(\cdot,\tau)|_q d\tau$.
The proof of Theorem \ref{big} is then concluded after we apply Gronwall's
lemma to (\ref{final}).

\vspace{.15in}
Now We turn our attention to the regularity issue of 
weak solutions of the IVP (\ref{dqg}) with $\alpha$ equal to
the critical index $\frac12$. The purpose of the next several theorems
is to show how far one can go toward
a regularity proof and what assumptions are needed to fill the gap.
 
\begin{thm}
Let $\alpha=\frac12$ and $s>0$. Assume that $\theta_0\in 
H^s\cap L^\infty$ and 
$f\in L^1([0,T];L^\infty)\cap L^2([0,T]; H^{s-\frac12})$. 
Consider a weak solution $\theta$ of the IVP (\ref{dqg}) and assume that 
\begin{equation}\label{assu}
|\theta(\cdot,\tau)|_\infty + |u(\cdot,t)|_\infty < \frac{\kappa}{C_0},
\end{equation}
where $C_0$ is a constant as in (\ref{key}).
Then $\theta$ is regular in the sense that 
\begin{equation}\label{reg}
\theta\in L^\infty([0,T]; H^s) \cap L^2([0,T];H^{s+\frac12})
\end{equation}
In particular, if $\theta_0$, $f$ and $u$ satisfy
\begin{equation}\label{assu1}
|\theta_0|_\infty + \int_0^t|f(\cdot,\tau)|_\infty d\tau
\le \frac{\kappa}{2C_0}\quad \mbox{and}\quad 
|u(\cdot,t)|_\infty < \frac{\kappa}{2C_0},
\end{equation}
then (\ref{reg}) holds.
\end{thm}
{\it Proof}.\quad The situation is different when $\alpha=\frac12$. The 
differential inequality (\ref{key}) 
only holds for $\beta=\alpha=\frac12$ and $q=\infty$. We have after replacing
$\beta$ with $\frac12$ and $q$ with $\infty$
$$
\frac{d}{dt} |\Lambda^s \theta|^2_2 + \kappa |\Lambda^{s+\frac12}
\theta|^2_2  \le \frac{1}{\kappa}|\Lambda^{s-\frac12}f|^2_2 
+ C_0(|\theta|_\infty + |u|_\infty) |\Lambda^{s+\frac12}\theta|_2^2
$$
Using the assumption (\ref{assu}), we have 
$$
\frac{d}{dt} |\Lambda^s \theta|^2_2 \le 
\frac{1}{\kappa}|\Lambda^{s-\frac12}f|^2_2,
$$
which gives (\ref{reg}). In view of the maximum principle
$$
|\theta(\cdot,t)|_\infty \le |\theta_0|_\infty 
+ \int_0^t |f(\cdot,\tau)|_\infty d\tau,
$$ 
it then follows that (\ref{assu1}) implies (\ref{assu}) and thus (\ref{reg}).

\vspace{.15in}
For the critical index $\alpha=\frac12$, the terms $|\theta|_\infty$ and 
$|u|_\infty$ have so far stood in the way of finding a regularity proof.
The problem of how to deal with $|\theta|_\infty$ and $|u|_\infty$ has to be
solved. In order to estimate $|\theta|_\infty$ and $|u|_\infty$, we first
prove the following lemma.
\begin{lemma}\label{FG}
Let $\Omega=[0,2\pi]^2$ and $F\in H^\sigma(\Omega)$ 
($\sigma>1$) be periodic. Then
\begin{equation}\label{four}
|F|_\infty \le C\left[1+ \|F\|_1 \,
\sqrt{log\left(1+\|F\|_\sigma^{\frac1{\sigma-1}}\right)}\right] 
\end{equation}
\end{lemma}
{\it Proof}.\quad Consider the Fourier transform $\widehat{F}$ of $F$.
For $R>0$, 
$$
|F(x)| \le \sum_{|k|\le R} |\widehat{F}(k)| + \sum_{|k\ge R}|\widehat{F}(k)|
\le \sum_{|k|\le R} |k|^{-1}|k| |\widehat{F}(k)| + 
\sum_{|k|\ge R} |k|^{-\sigma}|k|^{\sigma}|\widehat{F}(k)|
$$
$$
\le \left[\sum_{|k|\le R}|k|^{-2}\right]^\frac12 \, 
  \left[\sum_{|k|\le R} |k|^2|\widehat{F}(k)|^2\right]^\frac12
+ \left[\sum_{|k|\ge R}|k|^{-2\sigma}\right]^\frac12
 \left[\sum_{|k|\ge R} |k|^{2\sigma}|\widehat{F}(k)|^2\right]^\frac12
$$
$$
\le C \left[(\log(1+R))^\frac12 \|F\|_1 
+ \frac1{R^{\sigma-1}}\|F\|_\sigma\right] 
$$
In the above, we choose $R$ as 
$$
R\, = \|F\|_\sigma^{\frac1{\sigma-1}} 
$$
and we have after some manipulation the inequality (\ref{four}).

\vspace{.1in}
Now we can prove the following theorem, which can be reviewed as a ladder
theorem for the QG equation.
\begin{thm}
Let $\alpha=\frac12$, $T>0$ and $s>1$. Assume that $\theta_0\in H^s$ and 
$f\in L^2([0,T];H^{s-\frac12})$. Then any solution $\theta$ of the IVP 
(\ref{dqg}) satisfies for any $\sigma>1$
$$
\frac{d}{dt} |\Lambda^s \theta|^2_2 + \kappa |\Lambda^{s+\frac12}
\theta|^2_2  
$$
$$
\le \frac{1}{\kappa}|\Lambda^{s-\frac12}f|^2_2 
+ C \left[1+ \|\theta\|_1 \,
\sqrt{log\left(1+\|\theta\|_\sigma^{\frac1{\sigma-1}}\right)}\right]
|\Lambda^{s+\frac12}\theta|_2^2
$$ 
If we further assume that 
$\|\theta(\cdot,t)\|_\sigma$ is bounded and
small, say, for some constant $C$
$$
1+ \|\theta\|_1 \,
\sqrt{log\left(1+\|\theta\|_\sigma^{\frac1{\sigma-1}}\right)} \le C \kappa,
$$
then $\theta$ is regular in the sense that 
$$
\theta\in L^\infty([0,T];H^s) \cap L^2([0,T]; H^{s+\frac12}).
$$ 
\end{thm}

\vspace{.2in}
\section{Regularized QG Equation}
\setcounter{equation}{0}
\label{sec:3}

In this section we are concerned with the IVP for the regularized QG 
equation
\begin{equation}\label{rqg}
\left\{
\begin{array}{ll}
\theta_t  + u\cdot \nabla \theta 
+ \mu (-\Delta)^{\alpha} \theta_t = f,\quad 
& (x,t)\in \Omega\times [0,\infty),\\\cr
u=(u_1,u_2) = (-{\cal R}_2 \theta,\, {\cal R}_1\theta),\quad
& (x,t)\in \Omega\times [0,\infty),\\\cr
\theta(x,0) =\theta_0(x),\quad
& x\in \Omega,\cr
\end{array}
\right.
\end{equation}
where $0\le \alpha\le 1$ and $\mu>0$ are real numbers. The central issue 
is still whether or not weak solutions are regular for all 
time.  Because of the insignificant role of $f$ and for 
the sake of clarity of our presentation, we will set $f=0$ in the rest 
of this section.

\vspace{.15in}
The exposition below is organized as follows. We first construct a local 
solution in $H^s$ for $s>1$ and then produce explicit bounds on the norms 
of all derivatives. From this we infer that for $\alpha>\frac12$ 
the solutions are smooth and unique for all time. For $\alpha=\frac12$, 
global regularity is established under a weak assumption.
 
\vspace{.15in}
We reformulate the problem as an integral equation and then apply the 
Banach contraction mapping principle to prove local existence.
Now rewrite the first equation in (\ref{rqg}) in the form 
$$
(1+ \mu \Lambda^{2\alpha}) \,\theta_t = -(f- u\cdot\nabla\theta)
$$
and invert the operator $(1+\mu \Lambda^{2\alpha})$ subject to 
periodic boundary condition to obtain 
\begin{equation}\label{save}
\theta(x,t) =\theta_0(x) + \int_0^t G*(u\theta)(\tau) d\tau,
\end{equation}
where the convolution kernel $G(x)$ is defined through its Fourier 
transform  
\begin{equation}\label{kl}
{\widehat G} (k) ={\widehat G} ((k_1,k_2)) =\frac{i}{1+\mu |k|^{2\alpha}} 
\left[\begin{array}{c}k_1\\k_2\end{array}
\right].
\end{equation}

Now notice that for $\alpha\ge \frac12$ and any $s\ge 0$ 
\begin{equation}\label{est}
\begin{array}{c}
\|G* F\|_s =\left(\sum_k |k|^{2s} |\widehat{G*F}(k)|^2\right)^\frac12 
\\\cr
\le \left(\sum_{k} |\widehat{G}(k)|^2 |k|^{2s} 
|\widehat{F}(k)|^2\right)^\frac12 
\le \sup_{k}|\widehat{G}(k)|   \,\,\|F\|_s \le \frac{1}{\mu}\|F\|_s 
\end{array}
\end{equation}

\vspace{.1in}
We now state and prove a local existence result for smooth solutions 
of the IVP (\ref{rqg}).
\begin{thm}\label{local}
Let $\alpha\ge \frac12$ and assume that $\theta_0\in H^s$ for some $s>1$. 

\noindent (a)\,\,
There exists a $T=T(\|\theta_0\|_s)$ such that the IVP (\ref{rqg}) 
has a unique solution $\theta$ with $\theta\in L^\infty([0,T];H^s)$.

\noindent (b)\,\,
If $T_*$ is the supremum of the set of all $T>0$ such that (\ref{rqg})
has a solution in $L^\infty([0,T];H^s)$, then either $T_*=\infty$ or
$$
\limsup_{t\uparrow T_*}\|\theta(\cdot,t)\|_s =\infty.
$$
\end{thm}
{\it Proof}. \quad Write the equation (\ref{save}) symbolically as 
$\theta = A \theta$. $A$ is seen to be a mapping of the space 
$E=L^\infty([0,T];H^s)$ into itself, where $T>0$ is yet to be 
specified. 

\vspace{.1in}
Let $b=\|\theta_0\|_s$, and set $R=2b$. Define $B_R$ to be the ball with
radius $R$ centered at the origin in $E$. We now show that if $T$ is 
sufficiently small, then $A$ is a contraction map on $B_R$. 
Let $\theta$ and $\bar{\theta}$ be any two elements of $B_R$. Then we have
$$
\|A\theta-A\bar{\theta}\|_E =\left\|\int_0^t G*(u\theta-\bar{u}\bar{\theta}) 
d\tau\right\|_E 
$$
$$
\le T \|G*(u\theta-\bar{u}\bar{\theta})\|_E 
\le T \left[\|G*((u-\bar{u})\theta)\|_E \, 
+ \|G*(\bar{u}(\theta-\bar{\theta}))\|_E\right]
$$
Using (\ref{est}), we have 
\begin{equation}\label{diff}
\|A\theta-A\bar{\theta}\|_E \le \frac{T}{\mu} \sup_{t\in[0,T]}\left[
|\Lambda^s((u-\bar{u})\theta)|_2 + |\Lambda^s(\bar{u}(\theta-\bar{\theta}))|_2
\right]
\end{equation}
Then applying the calculus inequality (\ref{cal}), we find that 
\begin{equation}\label{ex1}
|\Lambda^s((u-\bar{u})\theta)|_2 
\le C\left(|\Lambda^s(u-\bar{u})|_2 |\theta|_\infty
+ |u-\bar{u}|_\infty |\Lambda^s\theta|_2 \right)
\end{equation}
and 
\begin{equation}\label{ex2}
|\Lambda^s(\bar{u}(\theta-\bar{\theta}))|_2 \le C\left( 
|\Lambda^s(\theta-\bar{\theta})|_2 |\bar{u}|_\infty
+ |\theta-\bar{\theta}|_\infty |\Lambda^s\bar{u}|_2\right)
\end{equation}
Since $s>1$, we have the Sobolev inequality 
\begin{equation}\label{sob}
|F|_\infty \le C \|F\|_s,\qquad\mbox{for some constant $C$}.
\end{equation}
Inserting (\ref{ex1}) and (\ref{ex2}) into (\ref{diff}) after 
applying (\ref{sob}), we obtain
\begin{equation}\label{almost}
\|A\theta-A\bar{\theta}\|_E \le 
\frac{T}{\mu} \sup_{t\in[0,T]}\left[\|\theta\|_s \|u-\bar{u}\|_s + 
\|\bar{u}\|_s \|\theta-\bar{\theta}\|_s\right] 
\end{equation}
Since $u$ and $\theta$, as well as $\bar{u}$ and $\bar{\theta}$, are 
related by the second equation in (\ref{rqg}), one has
$$
\|\bar{u}\|_s \le \|\bar{\theta}\|_s, \quad
\|u-\bar{u}\|_s \le \|\theta-\bar{\theta}\|_s 
$$
It then follows from (\ref{almost}) that 
$$
\|A\theta-A\bar{\theta}\|_E \le \frac{T}{\mu}(\|\theta\|_E + 
\|\bar{\theta}\|_E) \|\theta-\bar{\theta}\|_E \le \frac{2TR}{\mu}
\|\theta-\bar{\theta}\|_E.
$$
Also, since $A0=\theta_0$  
$$
\|A\theta\|_E =\|A\theta- A0 +\theta_0\|_E 
\le \|A\theta- A0\|_E + \|\theta_0\|_E 
\le \frac{2TR}{\mu}\|\theta\|_E + b
$$
Now choose 
$$T =\frac{\mu}{4R}\quad\mbox{so that} \quad \frac{2TR}{\mu}=\frac12.$$ 
It is indeed that $T$ depends only on $\|\theta_0\|_s$ (since $R=2
\|\theta_0\|_s$). For this choice of $T$, we have 
$$
\|A\theta-A\bar{\theta}\|_E \le \frac12 \|\theta-\bar{\theta}\|_E,\quad
\mbox{and}\quad \|A\theta\|_E \le R
$$
The conclusion of part (a) then follows from contraction
mapping principle.

\vspace{.1in}
To prove part (b), suppose on the contrary that $T_*<\infty$ and 
that there exists a number $M>0$ and a sequence $t_n$ approaching 
$T_*$ from below such that
$$
\|\theta(\cdot, t_n)\|_s \le M,\quad \mbox{for $n=1,2,3,\cdots$}
$$
By part (a), there exists some $T=T(M)$ such that the solution 
starting with any $\theta(x,t_n)$ is in $L^\infty([0,T];H^s)$.
Since $t_n$ approaches $T_*$, we can choose $t_n$ such that 
$t_n+T>T_*$. By extending $\theta$ to the interval $[0,t_n+T]$,
we obtain a solution of the IVP (\ref{rqg}) in $L^\infty([0,t_n+T];H^s)$.
But this contradicts the maximality of $T_*$.

\vspace{.15in}
Now we return to the central issue: can the solution obtained in Theorem
\ref{local} be extended for all time? This motivates us to explore the 
regularity properties of solutions. For $\alpha>\frac12$, it is indeed
the case that the solution we constructed in Theorem \ref{local} can
be extended for all time.  
\begin{thm}\label{high}
Let $\alpha>\frac12$ and assume that $\theta_0\in H^s$ for some $s>1$.
Then there exists a unique solution to the IVP (\ref{rqg}) which lies in
$L^\infty([0,\infty); H^s)$.
\end{thm}
{\it Proof}.\quad Let $\theta\in L^\infty([0,T];H^s)$ be the solution 
of the IVP (\ref{rqg}) we constructed in Theorem \ref{local}. Take the
inner product of $\Lambda^{2s-2\alpha}\theta $ with the first equation in 
(\ref{rqg}) 
$$
\frac12\frac{d}{dt} \left[|\Lambda^{s-\alpha}\theta|^2_2 + 
\mu\,|\Lambda^{s}\theta|^2_2 \right] 
= - (\Lambda^{2s-2\alpha}\theta, u\cdot \nabla \theta).
$$ 
The estimates for the term on the right hand side are similarly to those 
in the previous section, so we only give the most important lines here.
Instead of (\ref{key}), we have 
\begin{equation}\label{some}
\frac12\frac{d}{dt} \left[|\Lambda^{s-\alpha}\theta|^2_2 +
\mu\,|\Lambda^{s}\theta|^2_2 \right]
 \le C_0(|\theta|_q + |u|_q) |\Lambda^{s-\alpha+\beta}\theta|_2^2
\end{equation}
for any $2<q\le\infty$ and $\beta=1/2+1/q\le \alpha$.
Inserting (\ref{ugood}) and the following modified version of 
(\ref{bgood}) in Theorem \ref{big}
$$
|\Lambda^{s-\alpha+\beta}\theta|_2^2 \le 
\sqrt{\mu} |\Lambda^{s}\theta |_2^2 + \frac{C}{\sqrt{\mu}}
|\Lambda^{s-\alpha}\theta|^2_2
$$
into (\ref{some}), we obtain
$$
\frac12\frac{d}{dt} \left[|\Lambda^{s-\alpha}\theta|^2_2 +
\mu\,|\Lambda^{s}\theta|^2_2 \right]
 \le \frac{C}{\sqrt{\mu}}\left[\Lambda^{s-\alpha}\theta|^2_2 +
\mu\,|\Lambda^{s}\theta|^2_2 \right] 
$$ 
Gronwall's lemma then implies that for any $t>0$ 
$$
|\Lambda^{s-\alpha}\theta(\cdot,t)|^2_2 +
\mu\,|\Lambda^{s}\theta(\cdot,t)|^2_2 \,\le \left[
|\Lambda^{s-\alpha}\theta_0|^2_2 +\mu\,|\Lambda^{s}\theta_0|^2_2\right]\, 
e^{\frac{C}{\sqrt{\mu} }t}.
$$
This estimate indicates that the solution of the IVP (\ref{rqg}) with
$\alpha>\frac12$ remains bounded at later time. Therefore, by part (b)
of Theorem \ref{local}, $\theta\in L^\infty ([0,\infty); H^s)$.

\vspace{.15in}
Now we turn our attention to $\alpha=\frac12$. The following theorem 
asserts that if the solution $\theta$ loses its regularity at a later
time, then the maximum of $\theta$ or the maximum of $u$ necessarily
grows without a bound. 
\begin{thm}\label{hhh}
Let $\alpha=\frac12$ and $\theta$ be the solution constructed in 
Theorem \ref{local}. If 
$T_*$ is the supremum of the set of all $T>0$ such that $\theta$ in 
the class $L^\infty([0,T];H^s)$ and $T_*<\infty$, then one of 
the following  
$$
\int_0^{T_*} |\theta(\cdot,\tau)|_\infty d\tau =\infty,\qquad
\int_0^{T_*} |u(\cdot,\tau)|_\infty d\tau =\infty
$$ 
holds. In other words, $\theta$ can be extended beyond $T_*$ to $T_1$ if 
for some constant $M$ and all $T<T_1$  
$$
\int_0^T |\theta(\cdot,\tau)|_\infty d\tau <M,\quad\mbox{and}
\quad \int_0^T  |u(\cdot,\tau)|_\infty d\tau <M.
$$
\end{thm}
{\it Proof}.\quad The differential inequality (\ref{some}) is valid for all 
$\alpha>0$, therefore we have by choosing $\alpha=\frac12$, $q=\infty$ and 
$\beta=\frac12$
\begin{equation}\label{sss}
\frac12\frac{d}{dt} \left[|\Lambda^{s-\frac12}\theta|^2_2 +
\mu\,|\Lambda^{s}\theta|^2_2 \right]
\le C_0(|\theta|_\infty + |u|_\infty)|\Lambda^{s}\theta|_2^2,
\end{equation}
which implies 
$$
|\Lambda^{s-\frac12}\theta(\cdot,t)|^2_2 +
\mu\,|\Lambda^{s}\theta(\cdot,t)|^2_2 
$$
$$
\le \left[
|\Lambda^{s-\frac12}\theta_0|^2_2 +\mu\,|\Lambda^{s}\theta_0|^2_2\right]\,
\exp{\left[\frac{C}{\mu}\int_0^t [|\theta(\cdot,\tau)|_\infty 
+ |u(\cdot,\tau)|_\infty] d\tau\right]}.
$$
The conclusion of the theorem is then inferred from this inequality.

\vspace{.15in}
Because of the Sobolev inequality (\ref{sob}), the following result is an easy 
consequence of Theorem \ref{hhh}.
\begin{cor}
Let $\alpha=\frac12$ and $\theta\in L^\infty([0,T];H^s)$ ($s>1$)
 be the solution
constructed in Theorem \ref{local}. Assume that $\theta$ satisfies for $T_1>T$
$$
\int_0^{T_1} \|\theta(\cdot, \tau)\|_\rho d\tau < \infty,
$$
where $\rho>1$. Then we can extend $\theta$ to be a solution of the IVP 
(\ref{rqg}) in the class
$L^\infty([0,T_1];H^s)$.
\end{cor}

\vspace{.15in}
The following theorem concludes that no singularities in $\|\theta\|_s$ 
($s>1$) are possible before $\|\theta\|_1$ becomes unbounded.
\begin{thm}
Let $\alpha=\frac12$ and $\theta\in L^\infty([0,T];H^s)$ ($s>1$)
 be the solution
constructed in Theorem \ref{local}. For any $T_1>T$, if 
$\theta\in L^\infty([0,T_1];H^1)$, then $\theta$ can be extended to 
$L^\infty([0,T_1]; H^s)$.
\end{thm}  
{\it Proof}.\quad After applying Lemma \ref{FG} to control 
$|\theta|_\infty$ and $|u|_\infty$, we have from (\ref{sss}) that 
$$
\frac12\frac{d}{dt} \left[|\Lambda^{s-\frac12}\theta|^2_2 +
\mu\,|\Lambda^{s}\theta|^2_2 \right]
\le C\left[1+ \|\theta\|_1 \,
\sqrt{log\left(1+\|\theta\|_s^{\frac1{s-1}}\right)}\right]
|\Lambda^{s}\theta|_2^2
$$
If $\theta\in L^\infty([0,T_1];H^1)$, i.e., $\|\theta\|_1< C$, then the above
inequality becomes 
\begin{equation}\label{z}
\frac{d z }{dt} \le C z \,\sqrt{1+log(z)}  
\end{equation}
where we have set 
$$
z(t) =|\Lambda^{s-\frac12}\theta(\cdot,t)|^2_2 +
\mu\,|\Lambda^{s}\theta(\cdot,t)|^2_2.
$$
It follows from applying Gronwall's 
lemma to (\ref{z}) that $z\in L^\infty([0,T_1])$, which 
implies $\theta\in L^\infty([0,T_1];H^s)$.

\vspace{.2in}
\section{Weak solutions of the QG equation}
\setcounter{equation}{0}
\label{sec:4}

Onsager conjectured \cite{On} that weak solutions of the 3D Euler equations 
in a H\"{o}lder space $C^\gamma$ with exponent 
$\gamma>\frac13$ should conserve 
energy. In \cite{Ey0} Eyink proved energy conservation for 
weak solutions in a strong form of H\"{o}lder space 
$C^\gamma_*$ ($\gamma>\frac13$), in which the norm is defined in terms of 
absolute Fourier coefficients.  
Constantin, E and 
Titi \cite{CET} proved a sharp version of Onsager's conjecture in Besov space 
$B_3^{s,\infty}$ with $s>\frac13$.
In \cite{DR} Duchon and Robert explored possible sources for energy losses
and proved the existence of the so-called dissipative weak solutions 
for the 3D Euler equations. Their notion of dissipative weak solutions
can be regarded as a special type of 
dissipative solutions proposed by Lions \cite{Li}. For the 2D Euler equations,
DiPerna and Majda constructed weak solutions of the
vorticity-velocity formulation with data in $L^p(p>\frac43)$ and 
Eyink \cite{Ey2} validated such weak solutions as  
dissipative weak solutions in the sense of 
Duchon and Robert and argued for their relevance to the enstrophy 
cascade of 2D turbulence.  

\vspace{.15in}
The goal of this section is to prove Onsager's conjecture and extend
the notion of dissipative weak solutions to the IVP for the QG equation
\begin{equation}\label{qg}
\left\{
\begin{array}{ll}
\theta_t  + u\cdot \nabla \theta 
=0,\quad 
& (x,t)\in \Omega\times [0,\infty),\\\cr
u=(u_1,u_2) = (-{\cal R}_2 \theta,\, {\cal R}_1\theta),\quad
& (x,t)\in \Omega\times [0,\infty),\\\cr
\theta(x,0) =\theta_0(x),\quad
& x\in \Omega.\cr
\end{array}
\right.
\end{equation}

\vspace{.1in}
Understanding the zero-dissipation limits of the Navier-Stokes equations 
(\cite{CW1},\cite{CW2}) is 
crucial in hydrodynamics turbulence theory and we believe that similar
limits for the regularizations of the QG equation will be equally
important in the turbulence theory for quasi-geostrophic flows. At the
end of this section we show that the smooth solution of 
the regularized QG equation converges 
to that of the QG equation as $\mu\to 0$.

\vspace{.15in}
For $\theta_0\in L^2$ and any $T>0$, 
the IVP (\ref{qg}) has been shown to possess 
global weak solutions (in the distributional sense)
in $L^\infty([0,T];L^2)$ \cite{Re}. In the rest of this section
we will use $\phi$ to denote the standard mollifier in ${\Bbb R}^2$, 
$\phi^\epsilon(x)=\frac1{\epsilon^2}\phi\left(\frac{x}{\epsilon}\right)$
and $F^\epsilon=\phi^\epsilon*F$ for any tempered distribution $F$.

\vspace{.15in}
First we show that if the weak solution $\theta$ is also in the Besov space 
$B^{s,\infty}_3$ with $s>\frac13$, 
then $\theta$ conserves the $L^2$-norm. More details on Besov spaces can be 
found in \cite{St}, but here we only list some of the basic facts we will use. 
If $F\in B^{s,\infty}_p$ for $p\ge 1$, then
\begin{equation}\label{besov}
\begin{array}{l}
|F(\cdot+y)-F(\cdot)|_p \le C |y|^s \|F\|_{B^{s,\infty}_p};\\\cr
|F-F^\epsilon|_p \le C \epsilon^s \|F\|_{B^{s,\infty}_p};\\\cr
|\nabla F|_p \le C \epsilon^{s-1} \|F\|_{B^{s,\infty}_p}.
\end{array}
\end{equation}

\begin{thm}
Let $\theta\in L^\infty([0,T];L^2)$ be a weak solution of the IVP (\ref{qg})
corresponding to $\theta_0\in L^2$ and arbitrary $T>0$.
If further $\theta\in 
L^3([0,T];B^{s,\infty}_3)$ with $s>\frac13$, then for any $t\le T$ 
\begin{equation}\label{seek}
|\theta(\cdot,t)|_2 = |\theta_0|_2.
\end{equation}
\end{thm}
{\it Proof}.\quad The idea of the proof is similar to the one in \cite{CET}. 
If $\theta$ solves the IVP (\ref{qg}), then
$\theta^\epsilon$ satisfies the following equation
\begin{equation}\label{ooo}
\partial_t \theta^\epsilon + u^\epsilon\cdot\nabla \theta^\epsilon = 
\nabla \cdot \sigma^\epsilon,
\end{equation}
where $\sigma^\epsilon = u^\epsilon \theta^\epsilon\, -(u\theta)^\epsilon$.
It is easy to check that $\sigma^\epsilon$ can be represented by the 
formula 
$$
\sigma^\epsilon = (u-u^\epsilon)(\theta-\theta^\epsilon) - 
r^\epsilon(u,\theta),
$$
with 
$$
r^\epsilon(u,\theta) = \int \phi(y)[(u(x-\epsilon y)-u(x))(\theta(x-\epsilon y)
-\theta(x)] dy. 
$$
It then follows (\ref{ooo}) that 
\begin{equation}\label{cen}
\int_{\Omega} |\theta^\epsilon(\cdot,t)|^2 dx 
- \int_{\Omega}|\theta^\epsilon_0|^2 dx
=\int_0^t \int_\Omega \sigma^\epsilon \cdot \nabla \theta^\epsilon dx d\tau 
\end{equation}
We now estimate the term on the right hand side.
$$
\left|\int_0^t \int_\Omega \sigma^\epsilon \cdot 
\nabla \theta^\epsilon dx d\tau\right| 
\le C \int_0^t |\nabla \theta^\epsilon(\cdot,\tau)|_3 \cdot 
|\sigma^\epsilon(\cdot,\tau)|_{3/2} d\tau 
$$
$$
\le C \int_0^t |\nabla \theta^\epsilon(\cdot,\tau)|_3\cdot
\left[|u-u^\epsilon|_{3} |\theta-\theta^\epsilon|_3 + 
|r^\epsilon(u,\theta)|_{3/2}\right] d\tau
$$
Since $\theta$ and $u$ are related by the second equation in (\ref{qg}),
we have from (\ref{besov}) 
$$
|u-u^\epsilon|_3 \le |\theta-\theta^\epsilon|_3
\le C \epsilon^s \|\theta\|_{B^{s,\infty}_3}, \quad\mbox{and}
$$
$$
|r^\epsilon(u,\theta)|_{3/2} 
\le C \epsilon^{2s}\|\theta\|_{B^{s,\infty}_3}. 
$$
Therefore
$$
\left|\int_0^t \int_\Omega \sigma^\epsilon \cdot
\nabla \theta^\epsilon dx d\tau\right|
\le C \epsilon^{3s-1} \int_0^t \|\theta\|_{B^{s,\infty}_3}^3 d\tau
$$
and it approaches zero as $\epsilon\to 0$. The proof is then 
completed after letting $\epsilon\to 0$ in (\ref{cen}). 
 
\vspace{.15in}
The following theorem extends the notion of dissipative weak solutions
of the inviscid hydrodynamics equations (\cite{DR},\cite{Ey2}) 
to the QG equation.
\begin{thm}
Let $\theta\in L^\infty([0,T];L^2)$ be a weak solution of the IVP (\ref{qg})
corresponding to $\theta_0\in L^2$ and arbitrary $T>0$.
If a function 
$G: {\Bbb R} \to {\Bbb R}$ is $
C^2$,  strictly convex and has bounded derivative, then 
the equation 
\begin{equation}\label{noep}
\partial_t G(\theta) + u\cdot \nabla G(\theta) = - F(G,\theta),
\end{equation}
holds in the sense of distribution,
where $F(G,\theta)$ is the limit of  
$$G''(\theta_\epsilon) \nabla 
\theta_\epsilon\cdot ((u\theta)_\epsilon-u_\epsilon\theta_\epsilon) 
$$
in the sense of distribution.

\end{thm}
{\it Proof}.\quad Let $G\in C^2({\Bbb R},{\Bbb R})$. Multiplying
(\ref{ooo}) by $G'(\theta^\epsilon)$, we have
\begin{equation}\label{ep}
\partial_t G(\theta^\epsilon) + u^\epsilon\cdot\nabla G(\theta^\epsilon)
- \nabla\cdot (G'(\theta^\epsilon) \sigma^\epsilon) =-G''(\theta^\epsilon)
\sigma^\epsilon\cdot\nabla \theta^\epsilon.  
\end{equation}
Now we start showing that (\ref{ep}) converges to (\ref{noep}) in the 
distributional sense. Since $G$ has bounded derivative 
$$
|G'(x)|\,\le C 
$$
uniformly for all $x\in {\Bbb R}^2$, we have 
$$
|G(\theta^\epsilon)(\cdot,t) -G(\theta)(\cdot,t)|_2 
\le |G'|_\infty |\theta^\epsilon(\cdot,t)-\theta(\cdot,t)|_2 \quad \to 0 
$$
as $\epsilon\to 0$. Now we show that $u^\epsilon G(\theta^\epsilon)$ 
converges to $u G(\theta)$ in $L^1$, which can be deduced from the 
following estimate
$$
|u^\epsilon G(\theta^\epsilon)-u G(\theta)|_1 \le 
|(u^\epsilon-u)G(\theta^\epsilon)|_1  + |u(G(\theta^\epsilon)-G(\theta)|_1
$$
$$
\le |u^\epsilon-u|_2 |G(\theta^\epsilon)|_2 
+ |u|_2 |G(\theta^\epsilon)-G(\theta)|_2
$$
$$
\le C |\theta^\epsilon-\theta|_2 |G(\theta^\epsilon)|_2
+ |u|_2 |G'|_\infty |\theta^\epsilon-\theta|_2
$$
and the fact that $\theta^\epsilon\to \theta$ in $L^2$.  
It now remains to show that $G'(\theta^\epsilon) \sigma^\epsilon \to 0$ in
the distributional sense. 
$$
|G'(\theta^\epsilon) \sigma^\epsilon|_1 \le |G'|_\infty
|(u\theta)^\epsilon -u^\epsilon\theta^\epsilon|_1 
$$
$$
\le C |(u\theta)^\epsilon -u\theta|_1 + |u^\epsilon-u|_2|\theta|_2
+|u^\epsilon|_2|\theta^\epsilon-\theta|_2
$$
which approaches zero as $\epsilon\to 0$. Therefore the limit of the term
on the right hand side should also exist in the distributional sense.
This completes the proof of the theorem.

\vspace{.15in}
We now show that the two models (\ref{model1}) and (\ref{model3}) 
are close by examining the limit of equation (\ref{model3}) as 
$\mu\to 0$. First we recall that  
if $\theta_0\in H^s$ for $s\ge 3$, then 
the IVP (\ref{qg}) is known to have a unique smooth solution $\theta$ 
on a finite time interval satisfying $\theta\in L^\infty([0,T];H^s)$ 
{\cite{CMT} .   

\begin{thm}
Assume that $\{\theta_0^\mu\}_{\mu>0}$ and $\theta_0$ lie in $H^s$ with
$s\ge 3$. 
Then the difference $\theta^\mu-\theta$ between $\theta^\mu$ of the IVP 
(\ref{rqg}) with initial data $\theta_0^\mu$
and the solution $\theta$ of the IVP (\ref{qg})with initial data 
$\theta_0$ has the property 
$$
|\theta^\mu(\cdot,t)-\theta(\cdot,t)|_2^2 + \mu
\|\theta^\mu(\cdot,t)-\theta(\cdot,t)\|_\alpha^2 
$$
$$
\le \left(|\theta^\mu_0-\theta_0|_2^2 
+ \mu\|\theta^\mu_0-\theta_0\|_\alpha^2\right)
\exp\left[C\int_0^t (1+\|\theta(\cdot,\tau)\|_s)d\tau\right]
$$ 
$$
+ C \mu^2 \int_0^t \exp\left[C\int_0^{t-\zeta}
(1+\|\theta(\cdot,\tau)\|_s)d\tau\right]
\|\theta(\cdot,\zeta)\|_{2\alpha+1}^4 d\zeta.
$$
uniformly for $0\le t\le T$, where $C$ is a pure constant and $T$ is any fixed 
time less than the existence time for $\theta$.

\vspace{.1in}
In particular, if there is a constant $C$ such that
$$
|\theta^\mu_0-\theta_0|_2^2
+ \mu\|\theta^\mu_0-\theta_0\|_\alpha^2 
 \le C \mu^2 \quad\mbox{as $\mu\to 0$}, 
$$
then
$$
|\theta^\mu(\cdot,t)-\theta(\cdot,t)|_2^2 + \mu 
\|\theta^\mu(\cdot,t)-\theta(\cdot,t)\|_\alpha^2 \le C \mu^2
$$ 
uniformly for $0\le t\le T$, where $C$ is a constant depending only on 
$T$ and $\|\theta_0\|_s$.
\end{thm}
{\it Proof}.\quad The difference $w(x,t) = \theta^\mu(x,t) 
-\theta(x,t)$ solves the equation 
\begin{equation}\label{ddf}
w_t + u^\mu\cdot \nabla w + v\cdot \nabla\theta + 
\mu \Lambda^{2\alpha} w_t + \mu \Lambda^{2\alpha}\theta_t =0,
\end{equation}
where $v\,= u^\mu-u$.
Multiplying (\ref{ddf}) by $w$ and integrate over $\Omega$, we obtain 
$$
\frac12\frac{d}{dt} \int \left[w^2 + \mu |\Lambda^\alpha w|^2\right] =
-\int v\cdot \nabla\theta \cdot w  -\mu \int w \Lambda^{2\alpha}\theta_t,
$$ 
where the two terms on the right-hand side may be estimated as follows.
$$
-\int v\cdot \nabla\theta \cdot w \le |\nabla\theta(\cdot,t)|_\infty
|v|_2 \, |w|_2 
$$
Since $|\nabla\theta(\cdot,t)|_\infty \le C \|\theta(\cdot,t)\|_s$ and 
$|v|_2\le C |w|_2$, it follows that for $s>2$
$$
-\int v\cdot \nabla\theta \cdot w 
\le C \|\theta(\cdot,t)\|_s |w|_2^2.
$$ 
Noticing that $\theta_t=-u\cdot\nabla\theta$ and applying the 
calculus inequality (\ref{cal}), we can bound the second 
term by
$$
\begin{array}{ll}
\mu\left|\int w \Lambda^{2\alpha+1}(u\theta)\right|
& \le |w|_2^2  + \frac{\mu^2}{4} |\Lambda^{2\alpha+1}(u\theta)|^2_2\\\cr  
& \le |w|_2^2  + C \mu^2\left[|\Lambda^{2\alpha+1}u|_2|\theta|_\infty
+ |u|_\infty |\Lambda^{2\alpha+1}\theta|_2\right]^2\\\cr
& \le |w|_2^2  + C \mu^2 |\Lambda^{2\alpha+1}\theta|_2^2 \,\,
\|\theta\|_{2\alpha+1}^2 \\\cr
& \le |w|_2^2  + C \mu^2 \|\theta\|_{2\alpha+1}^4
\end{array}
$$ 
Collecting the above estimates, there appears
\begin{equation}\label{hoo}
\frac{d}{dt} \int \left[w^2 + \mu |\Lambda^\alpha w|^2\right]
\le C (1+\|\theta(\cdot,t)\|_s) \int w^2 
+ C \mu^2 \|\theta\|_{2\alpha+1}^4, 
\end{equation}
where the pure constant $C$ does not depend on $\mu$.
The desired result then follows from applying Gronwall's lemma to 
(\ref{hoo}).

\end{document}